\documentclass[a4,11]{article}

\makeatletter

\def\affil#1{\def\@affil{#1}}

\def\@maketitle{
\begin{center}
{\Large \@title \par}
\vspace{20pt}
{\normalsize \@author \par}
\vspace{15pt}
{\it\footnotesize \@affil \par}
\vspace{15pt}
\par\vskip 1.5em
}

\makeatother

\usepackage{amsmath,amsthm,amssymb,mathrsfs}
\usepackage[dvipdfmx]{graphicx}
 \usepackage[colorlinks=true]{hyperref}
\hypersetup{urlcolor=blue, citecolor=red}

\newtheorem{theorem}{Theorem}[section]
\newtheorem{lemma}[theorem]{Lemma} 
\newtheorem{proposition}[theorem]{Proposition} 
 
\newtheorem{remark}[theorem]{Remark}

\def\R{{\mathbb R}}

\def\Hb{{H^1(\Omega)}}
\def\Ho{{H^1_0(\Omega)}}
\def\Hi{{H^1_{\Gamma_i}(\Omega)}}
\def\Hfirst{{H^1_{\Gamma_1}(\Omega)}}
\def\Hsecond{{H^1_{\Gamma_2}(\Omega)}}
\def\Lo{{L^2(\Omega)}}

\def\dvg{\operatorname{div}}

\def\into{\int_\Omega}
\def\intg{\int_\Gamma}

\newcommand{\DS}{\displaystyle}

\title{\Large Boundary conditions for the Stokes problem and\\ a pressure-Poisson problem}

\author{\normalsize Kazunori Matsui}
\affil{Division of Mathematical and Physical Sciences,\\
Graduate School of Natural Science and Technology,\\ 
Kanazawa University, Kakuma, Kanazawa 920-1192 Japan}

\begin{document}
\maketitle
\footnotetext{
  {\it Email address:} {\tt first-lucky@stu.kanazawa-u.ac.jp}\\
  \indent~~{\it Date:} \today}

\begin{abstract}

  We consider a boundary value problem for the stationary Stokes problem
  and the corresponding pressure-Poisson equation.
  We propose a new formulation for the pressure-Poisson problem 
  with an appropriate additional boundary condition.
  We establish error estimates between solutions to the Stokes problem 
  and the pressure-Poisson problem in terms of the additional boundary condition.
  As boundary conditions for the Stokes problem, 
  we use a traction boundary condition and a pressure boundary condition 
  introduced in C. Conca et al (1994).

\end{abstract}

\section{Introduction}

Let $\Omega$ be a bounded connected open set of $\R^3$
with Lipschitz continuous boundary $\Gamma$.
We assume that there exist two relatively open subsets
$\Gamma_1$ and $\Gamma_2$ of $\Gamma$ such that
\[
  |\Gamma\setminus(\Gamma_1\cup\Gamma_2)|=0,\quad
  |\Gamma_1|,|\Gamma_2|> 0,\quad
  \Gamma_1\cap\Gamma_2=\emptyset,\quad
  \mathring{\overline{\Gamma_1}}=\Gamma_1,\quad
  \mathring{\overline{\Gamma_2}}=\Gamma_2,
\]
where $\overline{A}$ is the closure of $A\subset\Gamma$ with respect to $\Gamma$,
$\mathring{A}$ is the interior of $A$ with respect to $\Gamma$
and $|A|$ is the two dimensional Hausdorff measure.

The strong form of the Stokes problem is given as follows. 
Find $u^S:\Omega\rightarrow\R^3$ and $p^S:\Omega\rightarrow\R$ such that
\begin{align}\label{eq_oldS}\tag{S}
  \left\{\begin{array}{rcll}
  -\Delta u^S+\nabla p^S &=& F & \mbox{in }\Omega, \\
  \dvg u^S &=& 0 & \mbox{in }\Omega, \\
  u^S&=&0&\mbox{on }\Gamma_1,\\
  T_\nu(u^S,p^S)&=&t^b&\mbox{on }\Gamma_2,
  \end{array}\right.
\end{align}
holds, where $t^b:\Gamma_2\rightarrow\R^3$,
$\nu$ is the unit outward normal vector for $\Gamma$ and
\begin{align*}
  T_\nu(u,p)_i:=
  \sum^3_{j=1}\left(\frac{\partial u_i}{\partial x_j}+\frac{\partial u_j}{\partial x_i}\right)\nu_j-p\nu_i
  \quad\mbox{for all }i=1,2,3.
\end{align*}
The functions $u^S$ and $p^S$ are the velocity and the pressure of the flow
governed by (\ref{eq_oldS}), respectively.
We refer to \cite{Girault} and \cite{Temam} for the details on the Stokes problem
(i.e., physical background and corresponding mathematical analysis).
Taking the divergence of the first equation, we obtain
\begin{align}\label{ppstrong}
\dvg F=\dvg(-\Delta u^S+\nabla p^S)=-\Delta (\dvg u^S)+\Delta p^S=\Delta p^S.
\end{align}
This equation is often called the pressure-Poisson equation and 
is used in numerical schemes such as MAC (marker and sell), 
SMAC (simplified MAC) or the projection method 
(see, e.g., \cite{Amsden_Harlow,Chorin68,Cummins_Rudman,Guermond,mac1,Kim_Moin,mac2,Perot}).

We need an additional boundary condition for solving the equation (\ref{ppstrong}). 
In the real-would applications, the additional boundary condition
is usually given by using experimental or plausible values.
We consider the following boundary value problem for 
the pressure-Poisson equation:
Find $u^{PP}:\Omega\rightarrow\R^3$ and $p^{PP}:\Omega\rightarrow\R$ satisfying
\begin{align}\label{eq_PP}\tag{PP}
\left\{\begin{array}{rcll}
-\Delta u^{PP}-\nabla(\dvg u^{PP})+\nabla p^{PP} &=& F & \mbox{in }\Omega, \\
-\Delta p^{PP} &=& -\dvg F& \mbox{in }\Omega, \\
u^{PP}&=&0& \mbox{on }\Gamma_1,\\
{\DS\frac{\partial p^{PP}}{\partial\nu}}&=&g^b& \mbox{on }\Gamma_1,\\[8pt]
T_\nu(u^{PP},p^{PP})&=&t^b& \mbox{on }\Gamma_2,\\[8pt]
p^{PP}&=&p^b&\mbox{on }\Gamma_2,
\end{array}\right.
\end{align}
where $g^b:\Gamma_1\rightarrow\R$, $p^b:\Gamma_2\rightarrow\R$
are the data for the additional boundary conditions.
We call this problem the pressure-Poisson problem.
The second term $-\nabla(\dvg u^{PP})$ in the first equation of (\ref{eq_PP}) 
is necessary to treat the traction boundary condition 
in a weak formulation.
The idea of using (\ref{ppstrong}) instead of $\dvg u^S=0$ is useful
for calculating the pressure numerically in the Navier--Stokes problem.
For example, this idea is used in the MAC, SMAC and projection methods.

\begin{figure}[htp]
\begin{center}
  \includegraphics[width=0.5\linewidth,bb=-0 -0 335 145]{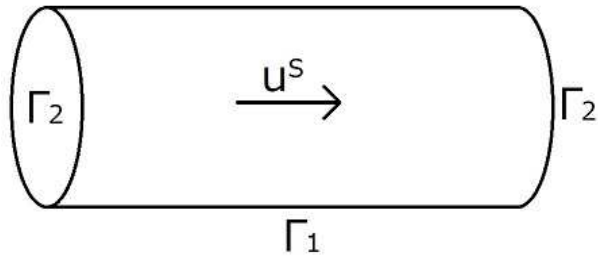}
  \caption{Image of a flow in a pipe}\label{flow}
  \end{center}
\end{figure}

In this paper, we establish error estimates between solutions 
for (\ref{eq_oldS}) and (\ref{eq_PP}) in terms of 
the additional boundary conditions.
As the boundary condition for the Stokes problem, 
we also consider the boundary condition introduced in \cite{Conca_etc94};
\begin{align}\label{BC_Conca}
\left\{\begin{array}{rcll}
u&=&0&\mbox{on }\Gamma_1,\\
u\times\nu&=&0&\mbox{on }\Gamma_2,\\
p&=&p^b&\mbox{on }\Gamma_2,
\end{array}\right.\end{align}
where ``$\times$'' is the cross product in $\R^3$
(see also \cite{Bertoluzza,Conca_etc95,Marusic}).
Since boundary conditions which contain a Dirichlet boundary condition
for the pressure often appear in engineering problems,
a comparison between (\ref{eq_PP}) and the Stokes problem with (\ref{BC_Conca}) 
is important.
For example, an end of pipe such as blood vessels or pipelines 
corresponds to the boundary $\Gamma_2$ (Fig.~\ref{flow}).

The organization of this paper is as follows.
In Section \ref{preliminaries} we introduce notations and symbols used in this work
and the weak form of these problems. We also prove the well-posedness of 
the problems (\ref{eq_oldS}) and (\ref{eq_PP}) and show several properties of them.
In Section \ref{sec_main1} we establish error estimates between solutions to 
the problems (\ref{eq_oldS}) and (\ref{eq_PP}) in terms of 
the additional boundary conditions.
Section \ref{sec_main2} is devoted to the study of the Stokes problem with 
the boundary condition (\ref{BC_Conca}).
We conclude this paper with several comments on future works 
in Section \ref{conclusion}.

%
\section{Preliminaries}\label{preliminaries}
\subsection{Notation}\label{notations}
%

We use the usual Lebesgue space $\Lo$ and Sobolev spaces 
$H^r(\Omega)$ $=W^{r,2}(\Omega)$ for a non-negative integer $r$,
together with their standard norms.
For spaces of vector-valued functions, we write $\Lo^3$, and so on.
The space $\Ho$ denotes the closure of $C^\infty_0(\Omega)$ in $\Hb$.
$\mathcal{D}'(\Omega)$ denotes the space of distributions on $\Omega$.
We set
\[\begin{array}{rcl}
  \Hi&:=&\{\psi\in\Hb\,|\,\psi=0\mbox{ on }\Gamma_i\}\quad(i=1,2),\\
  H&:=&{\DS\left\{\varphi\in\Hb^3\,\middle|\,\varphi=0\mbox{ on }\Gamma_1,
  \varphi\times\nu=0\mbox{ on }\Gamma_2\right\}.}
\end{array}\]

We also use the Lebesgue space $L^2(\Gamma)$ and Sobolev space $H^{1/2}(\Gamma)$
defined on $\Gamma$.
The norm $\|\eta\|_{H^{1/2}(\Gamma)}$ is defined by
\[
    \|\eta\|_{H^{1/2}(\Gamma)}
    :=\left(\|\eta\|^2_{L^2(\Gamma)}+\int_\Gamma \int_\Gamma
    \frac{|\eta(x)-\eta(y)|^2}{|x-y|^3}ds(x)ds(y)\right)^{1/2},
\]
where $ds$ denotes the surface measure of $\Gamma$.
For function spaces defined on $\Gamma_i~(i=1,2)$, we write $L^2(\Gamma_i)$, and so on.

We further set 
\[\begin{array}{c}
  {\DS 
  p_{,i}:=\frac{\partial p}{\partial x_i},\quad
  p_{,ij}:=\frac{\partial p}{\partial x_i\partial x_j},\quad
  D_{ij}(u):=\frac{1}{2}\left(u_{i,j}+u_{j,i}\right),}\quad
  {\DS D(u):D(\varphi):=\sum^3_{l,m=1}D_{lm}(u)D_{lm}(\varphi)}
\end{array}\]
for all $p:\Omega\rightarrow\R,u=(u_1,u_2,u_3):\Omega\rightarrow\R^3,
\varphi=(\varphi_1,\varphi_2,\varphi_3):\Omega\rightarrow\R^3$ and $i,j=1,2,3$.

%
\subsection{Preliminary results}\label{pre_results}
%

Let $\gamma_0\in {\mathscr L}(\Hb, H^{1/2}(\Gamma))$ be
the standard trace operator.
The trace operator $\gamma_0$ is surjective and satisfies 
$\mbox{Ker}(\gamma_0)=\Ho$ \cite[Theorem 1.5]{Girault}.
Let $\nu$ be the unit outward normal for $\Gamma$.
Since $\nu$ is a unit vector, 
$\Hb^3\ni u\mapsto u\cdot\nu:=(\gamma_0 u)\cdot\nu\in H^{1/2}(\Gamma)$
is a linear continuous map.
For all $u\in \Hb^3$ and $p\in\Hb$, the following Gauss divergence formula holds:
\[
\into u\cdot\nabla p+\into(\dvg u)p=\int_\Gamma(u\cdot\nu)p.
\]

For $i=1,2$, composition of the trace operator $\gamma_0$ and 
the restriction $H^{1/2}(\Gamma)\rightarrow H^{1/2}(\Gamma_i)$ 
is denoted by $\psi\mapsto\psi|_{\Gamma_i}$. 
This map is continuous from $\Hb$ to $H^{1/2}(\Gamma_i)$.
Since the kernel of this map is $\Hi$, there exists a constant $c>0$ such that
\[
    \|\psi\|_{\Hb/\Hi}\le c\|\psi|_{\Gamma_i}\|_{H^{1/2}(\Gamma_i)},
\]
where $\|\psi\|_{\Hb/\Hi}:=\inf_{q\in\Hi}\|\psi+q\|_\Hb$.
We simply write $\psi$ instead of $\psi|_{\Gamma_i}$
when there is no ambiguity.
We denote by $\langle\cdot,\cdot\rangle_{\Gamma_i}$ the duality pairing between
$H^{-1/2}(\Gamma_i):=H^{1/2}(\Gamma_i)^*$ and $H^{1/2}(\Gamma_i)$.
We remark that $\eta^*\in L^2(\Gamma_i)$ can be identified with 
a element of $H^{-1/2}(\Gamma_i)$ by
\[
  \langle \eta^*,\eta\rangle_{\Gamma_i}:=\int_{\Gamma_i}\eta^*\eta
  \qquad\mbox{for all }\eta\in H^{1/2}(\Gamma_i).
\]

For $u\in\Hb^3$ and $p\in\Hb$ satisfying $\Delta u+\nabla(\dvg u)\in\Lo$
and $\Delta p\in\Lo$, we set 
\[\begin{array}{rcll}
  \langle D(u)\nu,\varphi\rangle_{\Gamma_2}
  &:=&{\DS \into\left(D(u):D(\varphi)
  +\frac{1}{2}(\Delta u+\nabla(\dvg u))\varphi\right)}
  &\mbox{for all }\varphi\in\Hfirst^3,\\[8pt]
  \hspace{5pt}{\DS \left\langle\frac{\partial p}{\partial\nu},\psi\right\rangle_{\Gamma_1}}
  &:=&{\DS \into\left(\nabla p\cdot\nabla\psi+(\Delta p)\psi\right)}
  &\mbox{for all }\psi\in\Hsecond.
\end{array}\]
We remark that $u\in H^2(\Omega)^3$ and $p\in H^2(\Omega)$ satisfy
\[
  \langle D(u)\nu,\varphi\rangle_{\Gamma_2}
  =\int_{\Gamma_2}\left(\sum^3_{i,j=1}D_{ij}(u)\varphi_i \nu_j\right),\quad
  \left\langle\frac{\partial p}{\partial\nu},\psi\right\rangle_{\Gamma_1}
  =\int_{\Gamma_1}\frac{\partial p}{\partial\nu}\psi
\]
for all $\varphi\in\Hfirst^3$ and $\psi\in\Hsecond$.
For $u\in\Hb^3$ and $p\in\Hb$ satisfying $\Delta u+\nabla(\dvg u)\in\Lo$, we set 
\[
  \langle T_\nu(u,p)\nu,\varphi\rangle_{\Gamma_2}
  :=2\langle D(u)\nu,\varphi\rangle_{\Gamma_2}-\int_{\Gamma_2}p\varphi\cdot\nu
  \qquad\mbox{for all }\varphi\in\Hfirst^3.
\]

We recall the following five theorems which are necessary for
the existence and the uniqueness of a solution to the Stokes problem.

\begin{theorem}\label{Stokes}
  {\rm \cite[Corollary 4.1]{Girault}}
Let $(X,\|\cdot\|_X)$ and $(Q,\|\cdot\|_Q)$ be two real Hilbert spaces.
Let $a:X\times X\rightarrow\R$ and $b:X\times Q\rightarrow\R$ be
bilinear and continuous maps and let $f\in X^*$.
If there exist two constants $\alpha>0$ and $\beta>0$ such that
\begin{align*}\begin{array}{rcll}
a(v,v)&\ge&\alpha\|v\|^2_X&\mbox{for all }v\in V,\\
{\DS \sup_{v\in X\setminus\{0\}}\frac{b(v,q)}{\|v\|_{X}}}&\ge&\beta\|q\|_Q&\mbox{for all }q\in Q,
\end{array}\end{align*}
where $V=\{v\in X~|~b(v,q)=0\mbox{ for all }q\in Q\}$,
then there exists a unique solution $(u,p)\in X\times Q$ to the following problem:
\begin{align*}\left\{\begin{array}{rcll}
a(u,v)+b(v,p)&=&f(v)&\mbox{for all }v\in X,\\
b(u,q)&=&0&\mbox{for all }q\in Q.
\end{array}\right.\end{align*}
\end{theorem}

\begin{theorem}\label{infsup}
  {\rm \cite[Lemma 3.4]{Marusic}}
  There exists a constant $c>0$ such that
  \begin{align*}
      \sup_{\varphi\in\Hfirst^3\setminus\{ 0\}}\frac{{\DS\into q\dvg\varphi}}{\|\varphi\|_{\Hb^3}}
      \ge c\|q\|_{\Lo}
  \end{align*}
  for all $q\in\Lo$.
\end{theorem}

The following theorem is called Korn's second inequality.

\begin{theorem}[Korn's second inequality]\label{Korn}
  {\rm \cite[Lemma 5.4.18]{Trangenstein13}}
  There exists a constant $c>0$ such that
  \[
      \| D(\varphi)\|_{\Lo^{3\times 3}}
      \ge c\|\varphi\|_{\Hb^3}.
  \]
  for all $\varphi\in\Hfirst^3$.
\end{theorem}

The following embedding theorem is called Poincare's inequality.

\begin{theorem}[Poincare's inequality]\label{Poincare}
  {\rm \cite[Lemma 3.1]{Girault}}
  There exists a constant $c>0$ such that 
  \[
    \|\nabla\varphi\|_{\Lo^3}\ge c\|\varphi\|_\Lo
  \] 
  for all $\varphi\in\Hi~(i=1\mbox{ or }2)$.
\end{theorem}

The following embedding theorem plays an important role in the proof of
the existence and the uniqueness of the solution to the Stokes problem
with the boundary condition (\ref{BC_Conca}).

\begin{theorem}\label{thm_conca}
  {\rm \cite[Lemma 1.4]{Conca_etc94}}
  If $\Omega$ or $\Gamma$ satisfy one of the following conditions;
  \[
    \Gamma\mbox{ is of class }C^{1,1}
  \]
  or
  \[
      \Omega\mbox{ is a convex polyhedron},
  \]
  then there exists a constant $c>0$ such that 
  \[
    \|\nabla\times v\|_{\Lo^3}\ge c\| v\|_{\Hb^3}
  \] 
  for all $v\in H$ satisfying $\dvg v=0$.
\end{theorem}

%
\subsection{Weak formulations of (\ref{eq_PP}) and (\ref{eq_oldS})}
\label{weakform}
%

We start by defining the weak solution to (\ref{eq_PP}).
Throughout of this paper, we suppose the following conditions;
\begin{align}\label{BC_stokes}
  t^b\in H^{-1/2}(\Gamma_2)^3,\qquad
  F\in\Lo^3,
\end{align}
\begin{align}\label{BC_pp}
  g^b\in H^{-1/2}(\Gamma_1),\quad 
  p^b\in\Hb,\quad
  \dvg F\in\Lo.
\end{align}

\begin{lemma}
  For $u\in H^2(\Omega)^3,p\in H^1(\Omega)$ 
  and $\varphi\in\Hfirst^3$, 
  \[
    \into(-\Delta u-\nabla(\dvg u)+\nabla p)\cdot\varphi
    =\frac{1}{2}\into D(u):D(\varphi)-\into p\dvg\varphi
    -\langle t,\varphi\rangle_{\Gamma_2}
  \]
  holds, where $t:=T_\nu(u,p)$.
\end{lemma}
\noindent $\textit{ Proof.}$
We compute
\[\begin{array}{l}
    \quad{\DS \into(-\Delta u-\nabla(\dvg u)+\nabla p)\cdot\varphi}\\[8pt]
    ={\DS -\into\sum^3_{i,j=1}(u_{i,jj}+u_{j,ij})\varphi_i
    +\into\sum^3_{i=1} p_{,i}\varphi_i}\\[8pt]
    ={\DS \sum^3_{i,j=1}\left\{\into (u_{i,j}+u_{j,i})\varphi_{i,j}
     -\intg (u_{i,j}+u_{j,i})\varphi_i \nu_j\right\}
     +\sum^3_{i=1}\left\{-\into p\varphi_{i,i}
     +\intg p\varphi_i\nu_i\right\}}\\[8pt]
    ={\DS\frac{1}{2}\into\sum^3_{i,j=1}(u_{i,j}+u_{j,i})(\varphi_{i,j}+\varphi_{j,i})
     -\into\sum^3_{i=1} p\varphi_{i,i}
     -\intg\sum^3_{i=1} T_\nu(u,p)_i\varphi_i }\\[8pt]
    ={\DS \frac{1}{2}\into D(u):D(\varphi)-\into p\dvg\varphi
     -\langle t,\varphi\rangle_{\Gamma_2}},
\end{array}\]
which completes the proof.
\qed

For the second equation of (\ref{eq_PP}), taking $\psi\in\Hsecond$,  
we obtain
\[\begin{array}{rl}
    {\DS -\into(\dvg F)\psi}
    &={\DS -\into(\Delta p^{PP})\psi}\\[8pt]
    &={\DS -\intg\frac{\partial p^{PP}}{\partial\nu}\psi
    +\into\nabla p^{PP}\cdot\nabla\psi}\\[8pt]
    &={\DS -\langle g^b,\psi\rangle_{\Gamma_1}+\into\nabla p^{PP}\cdot\nabla\psi.}
\end{array}\]
Therefore, the weak form of (\ref{eq_PP}) becomes as follows.
Find $u^{PP}\in\Hfirst^3$ and $p^{PP}\in\Hb$ such that
\begin{align}\label{eq_PP2}\tag{PP'}\left\{\begin{array}{rcll}
    {\DS\frac{1}{2}\into D(u^{PP}):D(\varphi)-\into p^{PP}\dvg\varphi}
    &=&{\DS\into F\cdot\varphi-\langle t^b,\varphi\rangle_{\Gamma_2}}
    &\mbox{for all }\varphi\in\Hfirst^3,\\[8pt]
    {\DS\into\nabla p^{PP}\cdot\nabla\psi}
    &=&{\DS-\into(\dvg F)\psi+\langle g^b,\psi\rangle_{\Gamma_1}}
    &\mbox{for all }\psi\in\Hsecond,\\[8pt]
    p^{PP}&=&p^b\hspace{95pt}\mbox{on }\Gamma_2.
\end{array}\right.\end{align}

\begin{remark}
  If $(u^{PP},p^{PP})\in\Hfirst^3\times\Hb$ satisfies 
  $u^{PP}\in H^2(\Omega)^3,p^{PP}\in H^1(\Omega)$ and (\ref{eq_PP2}),
  then we have
  \[\left\{\begin{array}{rcll}
    {\DS\into(-\Delta u^{PP}-\nabla(\dvg u^{PP})+\nabla p^{PP}-F)\cdot\varphi}
    &=&{\DS\langle T_\nu(u^{PP},p^{PP})-t^b,\varphi\rangle_{\Gamma_2}}
    &\mbox{for all }\varphi\in\Hfirst^3,\\[8pt]
    {\DS\into(-\Delta p^{PP}+\dvg F)\cdot\psi}
    &=&{\DS\left\langle-\frac{\partial p^{PP}}{\partial\nu}+g^b,\psi\right\rangle_{\Gamma_1}}
    &\mbox{for all }\psi\in\Hsecond.
  \end{array}\right.\]
  Therefore, $(u^{PP},p^{PP})$ satisfies (\ref{eq_PP}).
\end{remark}
\medskip

Next, we define the weak formulation of (\ref{eq_oldS}).
For all $\varphi\in\Hfirst^3$, we obtain from the first equation of (\ref{eq_oldS}),
\[\begin{array}{rl}
    {\DS \into F\cdot \varphi}
    =&{\DS \into(-\Delta u^S+\nabla p^S)\cdot \varphi}\\[8pt]
    =&{\DS \into(-\Delta u^S-\nabla(\dvg u^S)+\nabla p^S)\cdot \varphi}\\[8pt]
    =&{\DS \frac{1}{2}\into D(u^S):D(\varphi)-\into p^S\dvg \varphi
    -\langle t^b,\varphi\rangle_{\Gamma_2}}.
\end{array}\]
Using this expression,
the weak form of the Stokes problem becomes as follows:
Find $(u^{S1},p^{S1})\in\Hfirst^3\times\Lo$ such that
\begin{align}\label{eq_S1}\tag{S1}\left\{\begin{array}{rcll}
    {\DS\frac{1}{2}\into D(u^{S1}):D(\varphi)-\into p^{S1}\dvg \varphi}
    &=&{\DS\into F\cdot \varphi -\langle t^b,\varphi\rangle_{\Gamma_2}}
    &\mbox{for all }\varphi\in\Hfirst^3,\\[8pt]
    {\DS -\into\psi\dvg u^{S1}}&=&0&\mbox{for all }\psi\in\Lo.
\end{array}\right.\end{align}

\begin{remark}
  If $(u^{S1},p^{S1})\in\Hfirst^3\times\Lo$ satisfies 
  $u^{S1}\in H^2(\Omega)^3,p^{S1}\in H^1(\Omega)$ and (\ref{eq_S1}),
  then we have
  \[\left\{\begin{array}{rcll}
    {\DS\into(-\Delta u^{S1}+\nabla p^{S1}-F)\cdot\varphi}
    &=&{\DS\langle T_\nu(u^{S1},p^{S1})-t^b,\varphi\rangle_{\Gamma_2}}
    &\mbox{for all }\varphi\in\Hfirst^3,\\[8pt]
    {\DS\into\psi\dvg u^{S1}}&=&0&\mbox{for all }\psi\in\Lo.\\
  \end{array}\right.\]
  Therefore, $(u^{S1},p^{S1})$ satisfies (\ref{eq_oldS}).
\end{remark}

%
\subsection{Well-posedness of (\ref{eq_PP2}), (\ref{eq_S1})}
%

We show the well-posedness of the problems (\ref{eq_PP2}) and (\ref{eq_S1})
in Theorem \ref{pp_thm} and \ref{s1u_thm}.

\begin{theorem}\label{pp_thm}
Under the conditions (\ref{BC_stokes}) and (\ref{BC_pp}),
there exists a unique solution $(u^{PP},p^{PP})\in\Hfirst^3\times\Hb$ satisfying (\ref{eq_PP2}).
\end{theorem}

\noindent $\textit{ Proof.}$
From the second and third equations of (\ref{eq_PP2}),
by using the Lax--Milgram theorem and Theorem \ref{Poincare},
$p^{PP}\in\Hb$ is uniquely determined.
Then $u^{PP}\in\Hb^n$ is also uniquely determined 
from the first equation of (\ref{eq_PP2}) by the Lax--Milgram theorem,
where the coercivity is guaranteed from Theorem \ref{Korn}.
\qed

\begin{theorem}\label{s1u_thm}
Under the condition (\ref{BC_stokes}),
there exists a unique solution $(u^{S1},p^{S1})\in\Hfirst^3\times\Lo$ satisfying (\ref{eq_S1}).
\end{theorem}

\noindent $\textit{ Proof.}$
By Theorems \ref{Korn} and \ref{Poincare},
the continuous bilinear form $\Hfirst^3\times\Hfirst^3\ni(u,\varphi)
\mapsto\into D(u):D(\varphi)\in\R$ is coercive. 
By Theorems \ref{Stokes} and \ref{infsup}, there exists a unique solution 
$(u^{S1},p^{S1})\in\Hfirst^3\times\Lo$ satisfying (\ref{eq_S1}).
\qed

We prove the following property of the solution to (\ref{eq_S1}).

\begin{proposition}\label{s1_ws}
If the weak solution $(u^{S1},p^{S1})\in\Hfirst^3\times\Lo$ to (\ref{eq_S1})
satisfies $p^{S1}\in\Hb$ and $\Delta p^{S1}\in\Lo$, then we have
\[
  \into\nabla p^{S1}\cdot\nabla\psi
  =-\into(\dvg F)\psi+\left\langle\frac{\partial p^{S1}}{\partial\nu},\psi\right\rangle_{\Gamma_1}
\]
for all $\psi\in\Hsecond$.
\end{proposition}
  
\noindent $\textit{ Proof.}$
From the second equation of (\ref{eq_S1}) and $u^{S1}\in\Hb$, 
$\dvg u^{S1}=0$ holds in $\Lo$.
From the first equation of (\ref{eq_S1}), we obtain
\[
  -\Delta u^{S1}+\nabla p^{S1}
  =-\Delta u^{S1}-\nabla(\dvg u^{S1})+\nabla p^{S1}=F
  \quad\mbox{in }{\mathscr D}'(\Omega).
\]
Taking the divergence, we get
\[
  \dvg F=\dvg(-\Delta u^{S1}+\nabla p^{S1})=-\Delta (\dvg u^{S1})+\Delta p^{S1}=\Delta p^{S1}
  \quad\mbox{in }{\mathscr D}'(\Omega).
\]
By the assumptions $\Delta p^{S1}\in\Lo$ and $\dvg F\in\Lo$,
$\Delta p^{S1}=\dvg F$ holds in $\Lo$.
Multiplying $\psi\in\Hsecond$ and integrating over $\Omega$,
we get
\[
    -\into(\dvg F)\psi=-\into (\Delta p^{S1})\psi
    =\into\nabla p^{S1}\cdot\nabla\psi
    -\left\langle\frac{\partial p^{S1}}{\partial\nu},\psi\right\rangle_{\Gamma_1},
\]
which is the desired result.
\qed

\section{The traction boundary condition}\label{sec_main1}

The purpose of this paper is to give an estimate of the difference
between the solutions of the Stokes problem and the pressure-Poisson problem.
Roughly speaking, from (\ref{ppstrong}) and the second equation of (\ref{eq_PP}),
$\Delta(p^S-p^{PP})=0$ holds. Hence, we get
\[
  \| p^S-p^{PP}\|_\Hb\lesssim(\mbox{
    difference between }p^S\mbox{ and }p^{PP}\mbox{ on }\Gamma),
\]
where $A\lesssim B$ means that there exists a constant $c>0$,
independent of $A$ and $B$, such that
$A\le cB$.
From (\ref{eq_oldS}) and the second equation of (\ref{eq_PP}), we have
\[
  -\Delta(u^S-u^{PP})=-\nabla(p^S-p^{PP}).
\]
We obtain
\[
  \| u^S-u^{PP}\|_{\Hb^3}\lesssim
  \|\nabla(p^S-p^{PP})\|_{\Lo^3}+(\mbox{
    difference between }p^S\mbox{ and }p^{PP}\mbox{ on }\Gamma).
\]
Therefore, we have
\[\begin{array}{rl}
  &\| u^S-u^{PP}\|_{\Hb^3}+\| p^S-p^{PP}\|_\Hb\\
  \lesssim&(\mbox{
    difference between }(u^S,p^S)\mbox{ and }(u^{PP},p^{PP})\mbox{ on }\Gamma).
\end{array}\]
In other words, if we have a good prediction for the boundary data, 
then (\ref{eq_PP}) is good approximation for (\ref{eq_oldS}).

In this section, we prove these types of estimates for the weak solutions.
Let the solutions of (\ref{eq_PP2}) and (\ref{eq_S1}) be denoted by
$(u^{PP},p^{PP})$ and $(u^{S1}, p^{S1})$, respectively.
First, we establish a lemma.

\begin{lemma}\label{lem_mix}
  If $p\in\Hb$, $f\in\Lo$ and $g\in H^{-1/2}(\Gamma_1)$ satisfy 
  \begin{align}\label{lem_mix1}
    \into\nabla p\cdot\nabla\psi=\into f\psi+\langle g,\psi\rangle_{\Gamma_1}
    \quad\mbox{for all }\psi\in\Hsecond,
  \end{align}
  then there exists a constant $c>0$ such that
  \[
      \|p\|_\Hb\le c\left(\|f\|_{\Lo}
      +\| g\|_{H^{-1/2}(\Gamma_1)}+\| p\|_{H^{1/2}(\Gamma_2)}\right).
  \]
\end{lemma}

\noindent ${\textit Proof.}$
Let $p_0\in\Hb$ such that $p_0-p\in\Hsecond$. 
Putting $\psi:=p-p_0$ in (\ref{lem_mix1}), we have
\begin{align*}\begin{array}{rcl}
  {\DS \|\nabla(p-p_0)\|^2_{\Lo^3}}
  &=&{\DS \into\nabla (p-p_0)\cdot\nabla(p-p_0)}\\[8pt]
  &=&{\DS \into f(p-p_0)+\langle g,p-p_0\rangle_{\Gamma_2}
  -\into\nabla p_0\cdot\nabla(p-p_0)}\\[8pt]
  &\le&{\DS \| f\|_\Lo\| p-p_0\|_{\Lo}
  +\| g\|_{H^{-1/2}(\Gamma_1)}\| p-p_0\|_{H^{1/2}(\Gamma_1)}}\\
  &&~{\DS+\|\nabla p_0\|_{\Lo^3}\|\nabla(p-p_0)\|_{\Lo^3}}\\[4pt]
  &\le&{\DS (\| f\|_\Lo+c_1\| g\|_{H^{-1/2}(\Gamma_1)}
  +\| p_0\|_{\Hb})\| p-p_0\|_{\Hb}.}
\end{array}\end{align*}
By Theorem \ref{Poincare}, there exists a constant $c_2>0$ such that
\[
  c_2\| p-p_0\|^2_{\Hb}
  \le(\| f\|_\Lo+c_1\| g\|_{H^{-1/2}(\Gamma_1)}+\| p_0\|_{\Hb})\| p-p_0\|_{\Hb}.
\]
Hence,
\[
  \| p-p_0\|_{\Hb}
  \le c_3(\| f\|_\Lo+\| g\|_{H^{-1/2}(\Gamma_1)}+\| p_0\|_{\Hb}).
\]
Since $\| p\|_\Hb-\| p_0\|_\Hb\le\| p-p_0\|_{\Hb}$, we obtain
\begin{align}\label{mix_mid1}
  \| p\|_{\Hb}
  \le c_4(\| f\|_\Lo+\| g\|_{H^{-1/2}(\Gamma_1)}+\| p_0\|_{\Hb}).
\end{align}
For all $p_0\in\Hb$ satisfying $p_0-p\in\Hsecond$, (\ref{mix_mid1}) holds.
Therefore,
\begin{align*}\begin{array}{rcl}
  {\DS \| p\|_{\Hb}}
  &\le&{\DS c_4\left(\| f\|_\Lo+\| g\|_{H^{-1/2}(\Gamma_1)}
  +\inf_{q\in\Hsecond}\| p+q\|_{\Hb}\right)}\\[8pt]
  &=&{\DS c_4(\| f\|_\Lo+\| g\|_{H^{-1/2}(\Gamma_1)}
  +\| p\|_{\Hb/\Hsecond})}\\
  &\le&{\DS c_5(\| f\|_\Lo+\| g\|_{H^{-1/2}(\Gamma_1)}
  +\|  p\|_{H^{1/2}(\Gamma_2)}).}\\
\end{array}\end{align*}
\qed

Using Lemma \ref{s1_ws}, we prove the following theorem
which is the main result of this section.

\begin{theorem}\label{thm_s1-p}
  If $p^{S1}\in\Hb$ and $\Delta p^{S1}\in\Lo$, 
  there exists a constant $c>0$ such that
  \begin{align}\label{ineq_s1-p}
    \|u^{S1}-u^{PP}\|_{\Hb^3}+\| p^{S1}-p^{PP}\|_{\Hb}
    \le c\left(\left\| \frac{\partial p^{S1}}{\partial\nu}-g^b\right\|_{H^{-1/2}(\Gamma_1)}
      +\| p^{S1}-p^b\|_{H^{1/2}(\Gamma_2)}\right).
  \end{align}
\end{theorem}

\noindent ${\textit Proof.}$
Using Proposition \ref{s1_ws}, we obtain from (\ref{eq_S1}) and (\ref{eq_PP}),
\begin{align}\label{s1-p} \left\{\begin{array}{rcll}{\DS
\frac{1}{2}\into D(u^{S1}-u^{PP}):D(\varphi)}
&=&{\DS\into (p^{S1}-p^{PP})\dvg\varphi }
&{\rm \mbox{for all }} \varphi \in\Hfirst^3, \\[8pt]
{\DS
\into\nabla(p^{S1}-p^{PP})\cdot\nabla\psi}
&=&{\DS\left\langle\frac{\partial p^{S1}}{\partial\nu}-g^b,\psi\right\rangle_{\Gamma_1}}
&{\rm \mbox{for all }} \psi \in \Hsecond.
\end{array}\right.\end{align}
Putting $\varphi:=u^{S1}-u^{PP}\in\Hfirst^3$ in (\ref{s1-p}), we get
\[\begin{array}{rl}
{\DS\frac{1}{2}\| D(u^{S1}-u^{PP})\|^2_{\Lo^{n\times n}}}
&={\DS \into (p^{S1}-p^{PP})\dvg(u^{S1}-u^{PP})}\\[8pt]
&\le \|p^{S1}-p^{PP}\|_{\Lo}\|\dvg(u^{S1}-u^{PP})\|_{\Lo}\\
&\le \sqrt{3}\|p^{S1}-p^{PP}\|_{\Hb}\|u^{S1}-u^{PP}\|_{\Hb^3}.
\end{array}\]
From Theorem \ref{Korn}, 
\begin{align}\label{usp_psp}
  \|u^{S1}-u^{PP}\|_{\Hb^3}\le c_1\|p^{S1}-p^{PP}\|_{\Hb}
\end{align}
holds for a constant $c_1>0$.
By the second equation of (\ref{s1-p}) and Lemma \ref{lem_mix},
there exists a constant $c_2>0$ such that
\[\begin{array}{rcl}
  {\DS \| p^{S1}-p^{PP}\|_{\Hb}}
  &\le&{\DS c_2\left(\left\| \frac{\partial p^{S1}}{\partial\nu}-g^b\right\|_{H^{-1/2}(\Gamma_1)}
    +\| p^{S1}-p^{PP}\|_{H^{1/2}(\Gamma_2)}\right)}\\[8pt]
  &\le&{\DS c_2\left(\left\| \frac{\partial p^{S1}}{\partial\nu}-g^b\right\|_{H^{-1/2}(\Gamma_1)}
    +\| p^{S1}-p^b\|_{H^{1/2}(\Gamma_2)}\right).}
\end{array}\]
Therefore, it holds that
\[
  \|u^{S1}-u^{PP}\|_{\Hb^3}+\| p^{S1}-p^{PP}\|_{\Hb}
  \le c_3\left(\left\| \frac{\partial p^{S1}}{\partial\nu}-g^b\right\|_{H^{-1/2}(\Gamma_1)}
    +\| p^{S1}-p^b\|_{H^{1/2}(\Gamma_2)}\right),
\]
for a constant $c_3>0$.
\qed

\section{Boundary condition involving pressure}\label{sec_main2}

Let $p^b\in\Hb$.
We consider the Stokes problem with the boundary condition (\ref{BC_Conca}):
\begin{align}\label{eq_oldS2}
  \left\{\begin{array}{rcll}
  -\Delta u^S+\nabla p^S &=& F & \mbox{in }\Omega, \\
  \dvg u^S &=& 0 & \mbox{in }\Omega, \\
  u&=&0&\mbox{on }\Gamma_1,\\
  u\times\nu&=&0&\mbox{on }\Gamma_2,\\
  p&=&p^b&\mbox{on }\Gamma_2.
  \end{array}\right.
\end{align}
In this section, we evaluate the difference between the solutions to
(\ref{eq_PP}) and (\ref{eq_oldS2}) as in (\ref{ineq_s1-p}).
First, we define the weak formulation of (\ref{eq_oldS2})
and prove the existence and the uniqueness of the weak solution.
Next, we prove a proposition and a lemma as preparation 
for the proof of our main theorem: Theorem \ref{thm_s2-p}.

We define the weak formulation of (\ref{eq_oldS2}).
Multiplying the first equation of (\ref{eq_oldS2}) by $v\in H$,
integrating by parts in $\Omega$, and using the second equation of (\ref{eq_oldS2}),
we obtain
\[
  \into F\cdot v
  =\into(\nabla\times u^S)\cdot(\nabla\times v)
  -\into p^S\dvg v+\int_{\Gamma_2}p^b v\cdot\nu,
\]
where we have used the following lemma.
\begin{lemma}
  For $u\in H^2(\Omega)^3,p\in H^1(\Omega)$ and $v\in H$, there holds
  \[
    \into(-\Delta u+\nabla(\dvg u)+\nabla p)\cdot v
    =\into(\nabla\times u)\cdot(\nabla\times v)
    -\into p\dvg v+\int_{\Gamma_2} pv\cdot\nu.
  \]
\end{lemma}
\noindent $\textit{ Proof.}$
We compute
\[\begin{array}{rl}
  &{\DS\into(-\Delta u+\nabla(\dvg u)+\nabla p)\cdot v}\\[8pt]
  =&{\DS\into(\nabla\times(\nabla\times u)+\nabla p)\cdot v}\\[8pt]
  =&{\DS \into(\nabla\times u)\cdot(\nabla\times v)
      -\intg((\nabla\times u)\times\nu)\cdot v
      -\into p\dvg v+\intg pv\cdot\nu}\\[8pt]
  =&{\DS \into(\nabla\times u)\cdot(\nabla\times v)
      -\int_{\Gamma}(\nu\times v)\cdot(\nabla\times u)
      -\into p\dvg v+\int_{\Gamma_2} pv\cdot\nu}\\[8pt]
  =&{\DS \into(\nabla\times u)\cdot(\nabla\times v)
      -\into p\dvg v+\int_{\Gamma_2} pv\cdot\nu.}
\end{array}\]
\qed

The weak form of the Stokes problem (\ref{eq_oldS2}) becomes as follows: 
Find $(u^{S2},p^{S2})\in H\times\Lo$ such that
\begin{align}\label{eq_S2}\tag{S2}\left\{\begin{array}{rcll}
    {\DS\into(\nabla\times u^{S2})\cdot(\nabla\times v)-\into p^{S2}\dvg v}
    &=&{\DS\into F\cdot v-\int_{\Gamma_2}p^b v\cdot\nu}
    &\mbox{for all }v\in H,\\
    {\DS -\into\psi\dvg u^{S2}}&=&0
    &\mbox{for all }\psi\in\Lo.
\end{array}\right.\end{align}

\begin{remark}
  If $(u^{S2},p^{S2})\in H\times\Lo$ satisfies 
  $u^{S2}\in H^2(\Omega)^3,p^{S2}\in H^1(\Omega)$ and (\ref{eq_S2}),
  then we have
  \[\left\{\begin{array}{rcll}
    {\DS\into(-\Delta u^{S2}+\nabla p^{S2}-F)\cdot v}
    &=&{\DS\int_{\Gamma_2}(p^{S2}-p^b) v\cdot\nu}
    &\mbox{for all }v\in H,\\[8pt]
    {\DS\into\psi\dvg u^{S2}}&=&0&\mbox{for all }\psi\in\Lo.
  \end{array}\right.\]
  Therefore, $(u^{S2},p^{S2})$ satisfies (\ref{eq_oldS2}).
\end{remark}
\medskip

We establish the well-posedness of this problem (\ref{eq_S2}) 
in the following theorem.

\begin{theorem}\label{s2u_thm}
  {\rm \cite[Theorem 1.5]{Conca_etc94}}
  For $F\in\Lo^3$ and $p^b\in\Hb$,
  under the hypotheses of Theorem \ref{thm_conca},
  there exists a unique solution $(u^{S2},p^{S2})\in H\times\Lo$ to (\ref{eq_S2}).
\end{theorem}

\noindent $\textit{ Proof.}$
We set 
\[
  a(u,v):=\into(\nabla\times u)\cdot(\nabla\times v),\quad
  b(v,q):=-\into q\dvg v,\quad
  f(v):=\into F\cdot v-\int_{\Gamma_2}p^b v\cdot\nu
\]
for all $u,v\in H$ and $q\in\Lo$.
Clearly, $a$ and $b$ are continuous and bilinear forms and $f\in H^*$. 
By Theorem \ref{thm_conca}, $a$ is coercive on 
$\{ v\in H~|~b(v,q)=0\mbox{ for all }q\in\Lo\}=\{ v\in H~|~\dvg v=0\}$.
By Theorem \ref{infsup}, $b$ satisfies the assumption of Theorem \ref{Stokes}.
Therefore, there exists a unique solution $(u^{S2},p^{S2})\in H\times\Lo$ 
to (\ref{eq_S2}) by Theorem \ref{Stokes}.
\qed

From here on, let the solutions of (\ref{eq_PP2}) and (\ref{eq_S2}) be denoted by
$(u^{PP},p^{PP})$ and $(u^{S2}, p^{S2})$, respectively.
The solution $(u^{S2}, p^{S2})$ to (\ref{eq_S2}) satisfies the following property.

\begin{proposition}\label{s2_ws}
If $\Delta u^{S2}+\nabla(\dvg u^{S2})\in\Lo^3$, $p^{S2}\in\Hb$ 
and $\Delta p^{S2}\in\Lo$, then 
\[\left\{\begin{array}{rcll}
  {\DS\frac{1}{2}\into D(u^{S2}):D(\varphi)-\into p^{S2}\dvg\varphi}
  &=&{\DS\into F\cdot\varphi-\langle T_\nu(u^{S2},p^{S2}),\varphi\rangle_{\Gamma_2}}
  &\mbox{for all }\varphi\in\Hfirst^3,\\[8pt]
  {\DS\into\nabla p^{S2}\cdot\nabla\psi}
  &=&{\DS-\into(\dvg F)\psi+\left\langle\frac{\partial p^{S2}}{\partial\nu},\psi\right\rangle_{\Gamma_1}}
  &\mbox{for all }\psi\in\Hsecond,\\[8pt]
  p^{S2}&=&p^b&\mbox{on }\Gamma_2.
\end{array}\right.\]
\end{proposition}
  
\noindent $\textit{ Proof.}$
From the second equation of (\ref{eq_S2}) and $u^{S2}\in\Hb$, 
$\dvg u^{S2}=0$ holds in $\Lo$.
From the first equation of (\ref{eq_S2}), we obtain
\begin{align}\label{up_S2}
  -\Delta u^{S2}-\nabla(\dvg u^{S2})+\nabla p^{S2}
  =-\Delta u^{S2}+\nabla(\dvg u^{S2})+\nabla p^{S2}=F
\end{align}
in ${\mathscr D}'(\Omega)$. 
By the assumptions $\Delta u^{S2}+\nabla(\dvg u^{S2})\in\Lo^3$, 
$p^{S1}\in\Hb$ and $\dvg F\in\Lo$, equation (\ref{up_S2}) holds in $\Lo$.
Multiplying $\varphi\in\Hfirst$ and integrating over $\Omega$,
we get
\[\begin{array}{rcl}
  {\DS \into F\cdot\varphi}
  &=&{\DS \into(-\Delta u^{S2}-\nabla(\dvg u^{S2})+\nabla p^{S2})\cdot\varphi}\\[8pt]
  &=&{\DS \frac{1}{2}\into D(u^{S2}):D(\varphi)-\into p^{S2}\dvg\varphi
  +\langle T_\nu(u^{S2},p^{S2}),\varphi\rangle_{\Gamma_2}.}
\end{array}\]
Taking the divergence of (\ref{up_S2}), we have
\[
  \Delta p^{S2}=\dvg F \quad\mbox{in }{\mathscr D}'(\Omega).
\]
By the assumptions $\Delta p^{S2}\in\Lo$ and $\dvg F\in\Lo$,
$\Delta p^{S2}=\dvg F$ holds in $\Lo$.
Multiplying $\psi\in\Hsecond$ and integrating over $\Omega$,
we get
\[
    -\into(\dvg F)\psi=-\into (\Delta p^{S2})\psi
    =\into\nabla p^{S2}\cdot\nabla\psi
    -\left\langle\frac{\partial p^{S2}}{\partial\nu},\psi\right\rangle_{\Gamma_1}.
\]
Multiplying (\ref{up_S2}) by $v\in H$ and integrating over $\Omega$,
we get
\[\begin{array}{rcl}
  {\DS \into F\cdot v}
  &=&{\DS \into(-\Delta u^{S2}+\nabla(\dvg u^{S2})+\nabla p^{S2})\cdot v}\\[8pt]
  &=&{\DS \into(\nabla\times u^{S2})\cdot(\nabla\times v)
  -\into p^{S2}\dvg v+\int_{\Gamma_2} p^{S2} v\cdot\nu.}
\end{array}\]
By the first equation of (\ref{eq_S2}), it holds that
\[
  \int_{\Gamma_2}p^{S2} v\cdot\nu
  = -\into(\nabla\times u^{S2})\cdot(\nabla\times v)
  +\into p^{S2}\dvg v+\into F\cdot v
  =\int_{\Gamma_2}p^b v\cdot\nu
\]
for all $v\in H$. Hence, $p^{S2}=p^b$ holds in $H^{1/2}(\Gamma_2)$.
\qed

We establish a lemma.

\begin{lemma}\label{lem_velo}
  If $u\in\Hfirst^3$, $p\in\Lo$ and $t\in H^{-1/2}(\Gamma_2)$ satisfy 
  \begin{align}\label{lem_velo1}
    \frac{1}{2}\into D(u):D(\varphi)
    =\into p\dvg\varphi-\langle t,\varphi\rangle_{\Gamma_2}
    \quad\mbox{for all }\varphi\in\Hfirst,
  \end{align}
  then there exists a constant $c>0$ such that
  \[
    \|u\|_{\Hb^3}\le c(\|p\|_{\Lo}+\| t\|_{H^{-1/2}(\Gamma_2)}).
  \]
\end{lemma}

\noindent ${\textit Proof.}$
Putting $\varphi:=u$ in (\ref{lem_velo1}), we obtain
\begin{align*}\begin{array}{rcl}
  {\DS \frac{1}{2}\|D(u)\|^2_{\Lo^{3\times 3}}}
  &=&{\DS \into p\dvg u-\langle t,u\rangle_{\Gamma_2}}\\[8pt]
  &\le&{\DS \| p\|_\Lo\|\dvg u\|_{\Lo}
  +\| t\|_{H^{-1/2}(\Gamma_2)}\| u\|_{H^{1/2}(\Gamma_2)}}\\
  &\le&{\DS (\sqrt{3}\| p\|_\Lo+c_1\| t\|_{H^{-1/2}(\Gamma_2)})\| u\|_{\Hb^3},}
\end{array}\end{align*}
for a constant $c_1>0$.
By Theorem \ref{Korn}, there exists a constant $c_2>0$ such that
\[
  \frac{c_2}{2}\| u\|^2_{\Hb^3}
  \le(\sqrt{3}\| p\|_\Lo+c_1\| t\|_{H^{-1/2}(\Gamma_2)})\| u\|_{\Hb^3}.
\]
Hence, we obtain the result with $c=\frac{2}{c_2}\max\{\sqrt{3},c_1\}$.
\qed

The next theorem is the main result of this section.

\begin{theorem}\label{thm_s2-p}
  If $\Delta u^{S2}+\nabla(\dvg u^{S2})\in\Lo^3$, $p^{S2}\in\Hb$ 
  and $\Delta p^{S2}\in\Lo$, then there exists a constant $c>0$ such that
  \[\begin{array}{rcl}
    \|p^{S2}-p^{PP}\|_{\Hb}&\le&
    {\DS c\left\|\frac{\partial p^{S2}}{\partial\nu}-g^b\right\|_{H^{-1/2}(\Gamma_1)},}\\[8pt]
    \|u^{S2}-u^{PP}\|_{\Hb^3}&\le&
    {\DS c\left(
      \left\|\frac{\partial p^{S2}}{\partial\nu}-g^b\right\|_{H^{-1/2}(\Gamma_1)}
      +\|t^{S2}-t^b\|_{H^{-1/2}(\Gamma_2)^3}\right),}
  \end{array}\]
  where $t^{S2}=T_\nu(u^{S2},p^{S2})$.
\end{theorem}

\noindent ${\textit Proof.}$
Using Proposition \ref{s2_ws}, we obtain from (\ref{eq_S2}) and (\ref{eq_PP}),
\begin{align}\label{s2-p} \left\{\begin{array}{rcll}{\DS
\frac{1}{2}\into D(u^{S2}-u^{PP}):D(\varphi)}
&=&{\DS\into (p^{S2}-p^{PP})\dvg\varphi-\langle t^{S2}-t^b,\varphi\rangle_{\Gamma_2}}
&{\rm \mbox{for all }} \varphi\in\Hfirst^3, \\[8pt]
{\DS
\into\nabla(p^{S2}-p^{PP})\cdot\nabla\psi}
&=&{\DS\left\langle\frac{\partial p^{S2}}{\partial\nu}-g^b,\psi\right\rangle_{\Gamma_1}}
&{\rm \mbox{for all }} \psi\in\Hsecond,\\[8pt]
p^{S2}-p^{PP}&=&0&\mbox{on }\Gamma_2,
\end{array}\right.\end{align}
where $t^{S2}=T_\nu(u^{S2},p^{S2})$.
By the second equation of (\ref{s2-p}) and Lemma \ref{lem_mix},
there exists a constant $c_1>0$ such that
\[\begin{array}{rcl}
  {\DS \| p^{S2}-p^{PP}\|_{\Hb}}
  &\le&{\DS c_1\left(\left\| \frac{\partial p^{S2}}{\partial\nu}-g^b\right\|_{H^{-1/2}(\Gamma_1)}
    +\| p^{S2}-p^{PP}\|_{H^{1/2}(\Gamma_2)}\right)}\\[8pt]
  &\le&{\DS c_1\left\| \frac{\partial p^{S2}}{\partial\nu}-g^b\right\|_{H^{-1/2}(\Gamma_1)}.}
\end{array}\]
By the first equation of (\ref{s2-p}) and Lemma \ref{lem_velo},
\[\begin{array}{rcl}
  {\DS \| u^{S2}-u^{PP}\|_{\Hb^3}}
  &\le&{\DS c_2\left(\|p^{S2}-p^{PP}\|_\Lo+\|t^{S2}-t^b\|_{H^{-1/2}(\Gamma_2)^3}\right)}\\[8pt]
  &\le&{\DS c_2\left(\|p^{S2}-p^{PP}\|_\Hb+\|t^{S2}-t^b\|_{H^{-1/2}(\Gamma_2)^3}\right)}\\[8pt]
  &\le&{\DS c_2\left(c_1\left\| \frac{\partial p^{S2}}{\partial\nu}-g^b\right\|_{H^{-1/2}(\Gamma_1)}
    +\|t^{S2}-t^b\|_{H^{-1/2}(\Gamma_2)^3}\right)}.
\end{array}\]
\qed

\section{Conclusion and future works}\label{conclusion}

We have proposed a new formulation for the pressure-Poisson problem (\ref{eq_PP}).
We have established error estimates between the solutions 
to (\ref{eq_PP2}) and (\ref{eq_S1}) in Theorem \ref{thm_s1-p} and 
between the solutions to (\ref{eq_PP2}) and (\ref{eq_S2}) in Theorem \ref{thm_s2-p}.
Theorem \ref{thm_s1-p} and \ref{thm_s2-p} state that 
if we have a good prediction for the boundary data ($g^b$ or $p^b$),
then the pressure-Poisson problem is a good approximation for the Stokes problem.

For problem (\ref{eq_S2}), a finite element scheme is proposed in \cite{Bertoluzza}
(under the assumption that $\Gamma_2$ is flat).
On the other hand, in many practical problems,
the projection method is more popular due to its easiness in implementation.
Numerical comparison of (\ref{eq_PP2}) and (\ref{eq_S2}) is 
one of our interesting future works from those points of view.

As another extension of our research,
generalization of our results to the Navier--Stokes problem
is important but is still completely open.

\bibliography{estokes_DCDS}
\end{document}